\documentclass[12pt]{article}
\usepackage{amsmath,amssymb}
\usepackage{bm}
\newtheorem{tm}{Theorem}[section]
\newtheorem{lm}[tm]{Lemma}

\newtheorem{re}[tm]{Remark}
\newtheorem{df}[tm]{Definition}

\newtheorem{pr}[tm]{Proposition}
 
 \newenvironment{demo}[1]{\par\smallskip\par\begin{trivlist}
\item[]{\bf #1}\ }{\end{trivlist}\par\smallskip\par}
\newcommand{\Proof}{\begin{demo}{{\it Proof.\ }}}
\newcommand{\qed}{\end{demo}}
\newcommand{\toy}{\ \rule[0em]{0.5ex}{1.8ex}}
\newcommand{\QED}{\toy\end{demo}}
\newcommand{\la}{\langle}
\newcommand{\ra}{\rangle}
\newcommand{\nn}{\nonumber}
\newcommand{\III}{{\vert \kern-.10em \vert \kern-.10em \vert}}
\newcommand{\ve}{\varepsilon}

\makeatletter
 
 \@addtoreset{equation}{section}
\makeatother
 
\begin{document}
\setlength{\baselineskip}{15pt} 
%%%%%%%%%%%%%%%%%%%%%%%%%%%%%%%%%%%%%%%%%%%%%%%%%%%%%%%%%%%%%%%%%%%%%%
%
\bibliographystyle{plain}
\title{
Malliavin differentiability of solutions of rough differential equations
\footnote{ Revised on 6 Mar 2014
{\bf Mathematics Subject Classification}: 60H07, 60H99, 60G15.
{\bf Keywords}: rough path theory, Malliavin calculus, 
 Gaussian process, 
}
}
%%%%%%%%%%%%%%%%%%%%%%%%%%%%%%%%%%%%%%%%%%%%%%%%%%%%%%%%%%%%%%%%%%
%%%%%%%%%%%%%%%%%%%%%%%%%%%%%%%%%%%%%%%%%%%%%%%%%%%%%%%%%%%%%%%%%%%
\author{  Yuzuru INAHAMA
\footnote{ 
%e-mail: ~\tt{inahama@math.nagoya-u.ac.jp}
%}
%\\
Graduate School of Mathematics,   Nagoya University
%\\
Furocho,  Chikusa-ku, Nagoya,  464-8602, Japan.
\quad
%\\
E-mail:~\tt{inahama@math.nagoya-u.ac.jp}
}
}
%\date{ \today }
\date{   }
% \pagestyle{empty}
%
% Start !!!
%
\maketitle
% \thispagestyle{empty}
%%%%\hspace{-5.5mm}

%\begin{center}
%{\bf Abstract}
%\end{center}
%
%
\begin{abstract}
In this paper we study rough differential equations 
driven by Gaussian rough paths from the viewpoint of Malliavin calculus. 
Under mild assumptions 
on  coefficient vector fields and underlying Gaussian processes,
we prove that 
solutions at a fixed time is smooth in the sense of Malliavin calculus.
Examples of Gaussian processes
include fractional Brownian motion with Hurst parameter larger than $1/4$.
\end{abstract}

%\vspace{5mm}
%%%%%%%%%%%%%%%%%%%%%%%%%%%%%%%%%%%%%%%%%%%%%%%%%%%%%%%%%%%%%%%%%%%%%%
%
%
\section{Introduction and main result}
%%%%%%%%%%%%%%%%%%%%%%%%%%%%%%%%%%%%%%%%%%%%%%%%%%%%%%%%%%%%%%%%%%%%%%
%
%\section{Setting and statement of main result}

%Let $(w_t)_{ 0 \le t \le T} = (w_t^1, \ldots, w_t^d)_{ 0 \le t \le T}$
%be a centered, continuous, $d$-dimensional Gaussian process 
%with i.i.d. components which start at $0$ ($T>0$).
%We denote by $\mu$ and ${\cal H}$ its law and its Cameron-Martin space, respectively.
%Then, $({\cal W}, {\cal H}, \mu)$ 
%becomes an abstract Wiener space, where ${\cal W}$ is the closure  of ${\cal H}$
%with respect to the usual sup-norm in $C_0 ([0,T], {\mathbb R}^d)$,
%the space of ${\mathbb R}^d$-valued continuous paths that start at $0$.

%%%%%%%%\vspace{10mm}
%%%%%%%

Let $V_{i}: {\mathbb R}^e \to {\mathbb R}^e$ be a vector field on ${\mathbb R}^e$
with sufficient regularity ($0 \le i \le d$)
and 
let $G\Omega_p ({\mathbb R}^d)$ be the geometric rough path space over ${\mathbb R}^d$
with $p$-variation topology.
For ${\bf x} \in G\Omega_p ({\mathbb R}^d)$, 
 Young pairing $({\bf x}, \bm{\lambda}) \in G\Omega_p ({\mathbb R}^{d+1})$ is well-defined, 
where $\lambda_t =t$.
We consider the following ordinary differential equation in the rough path sense
(rough differential equation, RDE) driven by $({\bf x}, \bm{\lambda})$;
\begin{equation}\label{rde.def}
dy_t = \sum_{i=1}^d  V_i ( y_t) dx_t^i + V_0 ( y_t) dt
\qquad
\qquad
\mbox{with given \quad $y_0 \in {\mathbb R}^e$.}
\end{equation}
If $V_i$'s are of $C_b^{[p] +1}$ (or of ${\rm Lip}(\gamma)$ for some $\gamma >p$), 
then a unique solution 
${\bf z} \in G\Omega_p ({\mathbb R}^{d+1} \oplus {\mathbb R}^e)$
which satisfies $\pi_1 {\bf z}  = ({\bf x}, \bm{\lambda})$ exists,
where $\pi_1$ is the projection onto the first component (i.e. ${\mathbb R}^{d+1}$-component)
 in the rough path sense.
${\bf y} := \pi_2 {\bf z} \in G\Omega_p ({\mathbb R}^e)$
is also called a solution.
We usually write ${\bf z}  = ({\bf x},  \bm{\lambda}, {\bf y})$, 
although the symbol on the right hand side is slightly misleading. 
We write $y_t = y_0 +{\bf y}^1_{0,t} \in {\mathbb R}^e$.
By Lyons' continuity theorem, 
$$
{\bf x} \mapsto ({\bf x},  \bm{\lambda}) \mapsto {\bf z} \mapsto {\bf y}
$$
is locally Lipschitz continuous with respect to $p$-variation topology.
The map in the middle is called Lyons-It\^o map (associated with $V_i$'s).
So far everything is deterministic.

%%%%%%%%\vspace{10mm}
%%%%%%%

Now we introduce a stochastic process.
Let $(w_t)_{ 0 \le t \le T} = (w_t^1, \ldots, w_t^d)_{ 0 \le t \le T}$
be a centered, continuous, $d$-dimensional Gaussian process 
with i.i.d. components
which start at $0$ ($T>0$).
We denote by $\mu$ and ${\cal H}$ its law and its Cameron-Martin space, respectively.
Then, $({\cal W}, {\cal H}, \mu)$ 
becomes an abstract Wiener space, where ${\cal W}$ is the closure  of ${\cal H}$
with respect to the usual sup-norm
on $C_0 ([0,T], {\mathbb R}^d)$,
the space of ${\mathbb R}^d$-valued continuous paths that start at $0$.

If $w =(w_t)$ admits a lift ${\bf w}$ as a $G\Omega_p ({\mathbb R}^d)$-valued random variable,
then $y_t=y_t ({\bf w}, \bm{\lambda})$ 
is a ${\mathbb R}^e$-valued Wiener functional defined on $({\cal W}, {\cal H}, \mu)$.
We will study differentiability of this Wiener functional in the sense of Malliavin calculus.

%%%%%%%%\vspace{10mm}
%%%%%%%

Let $R(s,t) = {\mathbb E} [w^1_s w^1_t]$ be the covariance function.
%
%Our assumption on Gaussian process $w$ is as follows.
%See Section 15.1, Friz and Victoir \cite{fvbk} for the definition of 2D $\rho$-variation.
%
If we assume $R (s,t)$
is of 2D $\rho$-variation for some $\rho \in [1, 2)$,
then the natural lift of $w$ exist as  
 the limit in $G \Omega_p ({\mathbb R}^d)$ of piecewise linear approximations of $w$
for any $p \in (2\rho, 4)$.
It  is denoted by ${\bf w}$ and  called a Gaussian rough path.
(For the definition of 2D $\rho$-variation, see Section 15.1, Friz and Victoir \cite{fvbk}.)

Our assumption on the Gaussian process $w =(w_t)$ is as follows.
This is called complementary Young regularity in p. 449, \cite{fvbk}.
(For a general theory of Young translation and pairing, see Section 9.4, \cite{fvbk}.)
\\
\\
{\bf (H)}: 
$R (s,t)$
is of 2D $\rho$-variation for some $\rho \in [1, 2)$.
Moreover, there exist $p \in  (2\rho, 4)$ and $q \in [1,2)$ such that (1)~
$1/p +1/q >1$ 
and (2)~ ${\cal H}$ is continuously embeded in $C_0^{q -var} ([0,T], {\mathbb R}^d)$,
the space of ${\mathbb R}^d$-valued continuous paths
of finite $q$-variation that start at $0$.
\begin{re}
Assumption {\bf (H)} holds if one of the following conditions holds:
\\
\noindent 
{\rm (i)}~$R$ is  of finite 2D $\rho$-variation for some $\rho \in [1,3/2)$.
\\
\noindent 
{\rm (ii)}~$R$ is  of finite 2D mixed $(1,\rho)$-variation for some $\rho \in [1,2)$.
(See \cite{fggr}.
Fractional Brownian motion with Hurst parameter $H \in (1/4, 1/2]$
satisfies this condition.)
\end{re}

%Under this assumption, 
%the natural lift of $w$ exist as  
% the limit of piecewise linear approximations of $w$
% and is denoted by ${\bf w}$.
%It is called a Gaussian rough path.

%%%%%%%%\vspace{10mm}
%%%%%%%

Now we state our main theorem in this paper.
(The proof is given in Subsection \ref{subsec.proof}.)
%
%
%Note that fractional Brownian motion 
%with Hurst parameter $H \in (1/4, 1/2]$ satisfies the assumption 
%of the theorem.  
%
We say a function $f$ defined on a domain in a Euclidean space 
is of $C_b^n$ if it is of $C^n$ and 
$f$ and its derivatives $\nabla^j f ~(1 \le j \le n)$ are bounded.
($C_b^{\infty}$ is defined in a similar way.)
For a real separable Hilbert space ${\cal K}$,
${\mathbb D}_{r, k} ({\cal K})$ stands for the ${\cal K}$-valued Gaussian-Sobolev 
space in the sense of Malliavin calculus
with the integrability index $r$ and the differentiability index $k$.
We set ${\mathbb D}_{\infty} ({\cal K}) 
= \cap_{1 <r <\infty} \cap_{ k \ge 0} {\mathbb D}_{r, k} ({\cal K})$ as usual.

%%%%%%%%%%%%%%%%%%%%%%%%%%%%%%%

\begin{tm}\label{tm.main}
Assume {\bf (H)} and  that the vector fields 
$V_i~(0 \le i \le d)$ are of $C_b^{\infty}$. 
We consider RDE (\ref{rde.def}) with ${\bf x} ={\bf w}$.
Then, for any $t \in [0,T]$, 
$y_t \in {\mathbb D}_{\infty} ({\mathbb R}^e)$,
that is, 
$y_t$ is smooth in the sense of Malliavin calculus.
\end{tm}

%%%%%%%%%%%%%%%%%%%%%%%%%%%%%%%%%

Malliavin calculus for the solution of an RDE driven by a Gaussian rough path
 was started by Cass, Friz, and Victoir \cite{cfv, cf}.
Since then many papers have been written 
and it is now a very active research topic
(see \cite{boz1, boz2, chlt, dri, hp, ht} among others).
If one wants to analyze the solution by means of Malliavin calculus,
showing smoothness of $y_t$ is a crucial step.
There are of course preceding results for special cases. 
%
%
%%%%%%%%
%\vspace{10mm}
%%%%%%%
%\begin{re}
%{\rm (i)}~
%
The case of fractional Brownian rough path with 
Hurst parameter  $H \in (1/3, 1/2]$
was shown by Hairer-Pillai \cite{hp} with fractional calculus.
(There are some results when the coefficient vector fields satisfy a special 
Lie-bracket condition. See \cite{dri, ht}.)
To the author's knowledge, however, 
this problem was not solved in a sufficiently general form.

%{\rm (ii)}~ In \cite{chlt} the authors proved this problem for a rather general class of 
%Gaussian rough paths.
%There seems to be a small gap in their proof, however.   
%At the moment of writing, it is somewhat unclear whether the proof can be fixed 
%or to what extent it really works.
%
%
%\\
%{\rm (iii)}~
%
\begin{re}
When $H=1/2$, 
fractional Brownian motion is Brownian motion
in the usual sense and $(y_t)$ is equal to the solution of 
the corresponding stochastic differential equation of
Stratonovich-type.
Our main theorem specialized in this case is slightly weaker than the classical
result in Malliavin calculus, where the condition on $V_i$ is as follows:
"For all $n=1,2,\ldots$ and $0 \le i \le d$,
$ \|\nabla^n V_i \|$
is bounded."
($V_i$ itself is allowed to have linear growth.)
However, since Bailleul \cite{be} recently solved RDE with such coefficients,
it might be possible to extend Theorem \ref{tm.main} to include such a case 
with existing methods. 
\end{re}

%%%%%%%%\vspace{10mm}
%%%%%%%

To prove our main theorem, 
we only use basic results in Malliavin calculus and 
some deep results in rough path theory, in particular, Friz and Victoir \cite{fv} and  
Cass, Litterer, and Lyons \cite{cll}.
(We do not use fractional calculus.)
The former one proves 
existence and basic properties of Gaussian rough paths
and 
the latter one proves 
integrability of the Jacobian process, which appears in an explicit expression of $D^n y_t$.
Our assumption {\bf (H)} is basically the same as the one in  \cite{cll}.

%%%%%%%%\vspace{10mm}
%%%%%%%

The rest of this paper is devoted to proving Theorem \ref{tm.main} above.
From now on, we will assume the following: 
(i)~Without loss of generality, we may assume $T=1$ and $y_0 =0$.
(ii)~Since the drift term always behaves nicely, we assume $V_0 =0$ for the sake of simplicity.
Therefore, we will study the following RDE with a random driving noise
${\bf x} ={\bf w}$; 
\begin{equation}
\label{rde.def2}
%\nn
dy_t = \sum_{i=1}^d  V_i ( y_t) dx_t^i 
\qquad
\qquad
\mbox{with \quad $y_0 =0\in {\mathbb R}^e$.}
\end{equation}
In this case, $y_t = {\bf y}^1_{0,t}$
and the domain of Lyons-It\^o map is simply $G\Omega_p ({\mathbb R}^{d})$, 
not $G\Omega_p ({\mathbb R}^{d+1})$.

%%%%%%%%%%%%%%%%%%%%%%%%%%%%%%%%%%%%%%%%%%%%%%%%%%%%%%%%%%%%%%%%%%%%%%
%\newpage
%\section{Proof of  main theorem}
%
\section{Heuristics}

In this section 
we give a heurtistic argument on how to estimate Hilbert-Schmidt norms of 
the $n$th ${\cal H}$-derivative $D^n y_t$ for $n=1,2$, 
so that the reader could easily understand our strategy.
Note that the contents of this section are not mathematically rigorous.

Consider the following ODE (or RDE) driven by $w$;
\begin{equation}
\label{9heu_y.eq}
%\nn
dy_t = \sum_{i=1}^d  V_i ( y_t) dw_t^i 
\qquad
\qquad
\mbox{with \quad $y_0 =0\in {\mathbb R}^e$.}
\end{equation}
The Jacobian process and its inverse are given by 
\begin{eqnarray}
dJ_t &=& \sum_{i=1}^d  \nabla V_i ( y_t) \cdot J_t dw_t^i 
\qquad
\qquad
\mbox{with \quad $J_0 ={\rm Id}_e \in {\rm Mat}(e,e)$.}
\label{heu_J.eq}
\\
dK_t &=& - \sum_{i=1}^d K_t  \cdot \nabla V_i ( y_t)  dw_t^i 
\qquad
\qquad
\mbox{with \quad $K_0 ={\rm Id}_e \in {\rm Mat}(e,e)$.}
\label{heu_K.eq}
\end{eqnarray}
Here, 
"$\cdot$" stands for the matrix multiplication and
$\nabla V_i $ is regarded as ${\rm Mat}(e,e)$-valued.
In fact, $K_t =J_t^{-1}$.

If the system of ODEs (\ref{9heu_y.eq})--(\ref{heu_K.eq})
is interpreted in the rough path sense,
then $({\bf w}, {\bf y}, {\bf J}, {\bf K})$
is a well-defined random rough path on ${\mathbb R}^d \oplus {\mathbb R}^e \oplus {\rm Mat}(e,e)^{\oplus 2}$.
If ${\bf w}$ is not very bad, we can use Cass-Litterer-Lyons' integrability 
criterion for $J$ and $K$ to conclude that
(any component of) any level path of $({\bf w}, {\bf y}, {\bf J}, {\bf K})$
has moments of all order.

Next, let us see what derivatives of $y_t$ look like.
For brevity, we write $\sigma = [V_1, \ldots, V_d]$, which is $e \times d$ matrix.
(Then, (\ref{9heu_y.eq}) is simply $dy_t = \sigma ( y_t) dw_t$.)
By formal differentiation of (\ref{9heu_y.eq}) in the direction of $h \in {\cal H}$, we obtain
\begin{equation}
d D_h y_t 
= 
\nabla \sigma ( y_t) \la D_h y_t , dw_t \ra
+ 
\sigma ( y_t) dh_t
\qquad
\qquad
\mbox{with \quad $ D_h y_0 =0\in {\mathbb R}^e$.}
\label{heu_1.eq}
\end{equation}
By the method of variation of constants, (\ref{heu_1.eq}) is equivalent to
\begin{equation}
 D_h y_t 
= 
J_t \int_0^t K_s  \sigma ( y_s) dh_s.
\label{heu_1int.eq}
\end{equation}

In  a similar way, we have for $h, k \in {\cal H}$ that 
\begin{eqnarray}
d D^2_{h,k} y_t 
&=& 
\nabla \sigma ( y_t) \la D^2_{h,k} y_t , dw_t \ra
+ 
\nabla^2 \sigma ( y_t) \la D_h y_t, D_k y_t ,dw_t\ra
\nn\\
&&+
\nabla \sigma ( y_t) \la D_k y_t ,dh_t\ra 
+
\nabla \sigma ( y_t) \la D_h y_t ,dk_t\ra 
\quad
\mbox{with  $ D^2_{h,k} y_0 =0\in {\mathbb R}^e$.}
\label{heu_2.eq}
\end{eqnarray}
This is equivalent to
\begin{eqnarray}
  D^2_{h,k} y_t 
&=& 
J_t \int_0^t K_s \bigl\{ 
\nabla^2 \sigma ( y_s) \la D_h y_s, D_k y_s ,dw_s \ra
\nn\\
&&
\qquad\qquad
+
\nabla \sigma ( y_s) \la D_k y_s ,dh_s\ra 
+
\nabla \sigma ( y_s) \la D_h y_s ,dk_s\ra  
\bigr\}.
\label{heu_2int.eq}
\end{eqnarray}

Now let us estimate $Dy_t$.
Since $h$ is $q$-variational, $K_s  \sigma ( y_s)$ is $p$-variational,
and $1/p +1/q >1$,
the right hand side of (\ref{heu_1int.eq}) is Young integral.
So,  we get 
$$
| D_h y_t|
\lesssim \|J\|_{\infty}  \| K_{\cdot} \cdot \sigma ( y_{\cdot})\|_{p -var} \|h\|_{q-var}
\lesssim
\|J\|_{\infty}  \| K_{\cdot}\|_{p -var}  \|  y_{\cdot} \|_{p -var} \|h\|_{{\cal H}}.
$$
Hence, $\|Dy_t\|_{{\cal H}^*} \lesssim
\|J\|_{\infty}  \| K_{\cdot}\|_{p -var}  \|  y_{\cdot} \|_{p -var}$
and 
$\|Dy_t\|_{{\cal H}^*}$ has moments of all order.

The above proof may look simple and good.
However, 
if we continue to argue in this way, we are in trouble even when $n=2$.
It is of course  possible to prove  that
$$
| D^2_{h,k} y_t|
\lesssim 
\mbox{(a polynomial in $p$-variation norm of some Wiener functionals)}
\times \|h\|_{{\cal H}} \|k\|_{{\cal H}}
$$
in a similar way as above.
This is basically an estimate of the operator norm of $D^2 y_t$, however, 
and we cannot get an estimate of 
 Hilbert-Schmidt norm $\| D^2 y_t\|_{ {\cal H}^* \otimes {\cal H}^* }$  so easily from this.

%%%%%%%%\vspace{10mm}
%%%%%%%

To overcome this difficulty, we will "double the dimension."
Let 
$(b_t) = (b_t^1, \ldots, b_t^d)$
which has the same law as $(w_t)$'s and is independent of $(w_t)$. 
Let us consider $2d$-dimensional Gaussian process 
$(w_t; b_t)_{0 \le t \le 1}$ from now on.

The expectation with respect to $w$-variable and $b$-variable
are
denoted by ${\mathbb E}'$ and $\hat{\mathbb E}$, respectively.
The expectation with respect to $(w, b)$-variable 
is of course the product ${\mathbb E} = {\mathbb E}' \times \hat{\mathbb E}$.

We formally replace $h \in {\cal H}$ in $D^n_{h,h, \ldots,h} y_t$ 
with  $b$ and denote it by $\Xi_n (w,b)$.
For $n=1,2$,
we have from (\ref{heu_1int.eq}) and (\ref{heu_2int.eq}) 
that
\begin{eqnarray}
\Xi_1 (w,b)_t
&=&
J_t \int_0^t K_s  \sigma ( y_s) db_s.
\label{heu_xi1.eq}
\\
\Xi_2 (w,b)_t
&=&
J_t \int_0^t K_s \bigl\{ 
\nabla^2 \sigma ( y_s) \la   \Xi_1 (w, b)_s,   \Xi_1 (w, b)_s, dw_s \ra
\nn\\
&&
\qquad\qquad\qquad
+
2\nabla \sigma ( y_s) \la  \Xi_1 (w, b)_s ,db_s\ra  
\bigr\}.
\label{heu_xi2.eq}
\end{eqnarray}

We will fix $w$ for a while.
Then, $\Xi_n (w, \,\cdot\,)_t$ belongs to the inhomogeneous Wiener chaos of order $n$
for $n=1,2$
(i.e., $\Xi_n (w, b)_t$ is a "polynomial" of order $n$ in $b$-variable).
Then,  simple formal computations yield
\begin{eqnarray*}
\hat{D}_h \Xi_1 (w,b)_t
&=&
\lim_{\ve \to 0} \frac{\Xi_1 (w,b+ \ve h)_t - \Xi_1 (w,b)_t}{\ve}
=
D_h y_t (w),
\\
\hat{D}^2_{h,k} \Xi_2  (w,b)_t
&=&
2 D^2_{h,k} y_t (w).
\end{eqnarray*}
Here, $\hat{D}$ denotes the derivative with respect to $b$-variable.
Note that the right hand sides are both constant in $b$.

Recall that all ${\mathbb D}_{2,k} ({\mathbb R}^e)$-norms ($k=0,1,\ldots$)  are equivalent on 
each fixed inhomogeneous  Wiener chaos.
Hence, we have
\begin{eqnarray*}
\|D y_t (w) \|_{{\cal H}^* \otimes {\mathbb R}^e} 
&=& 
\hat{\mathbb E} [ \|\hat{D} \Xi_1 (w, \,\cdot\,)_t \|_{{\cal H}^* \otimes {\mathbb R}^e}^2]^{1/2}
\lesssim
\| \Xi_1 (w, \,\cdot\,)_t \|_{ {\mathbb D}_{2,1} }
\lesssim
\| \Xi_1 (w, \,\cdot\,)_t \|_{ L^{2} }
,
\\
\|D^2 y_t (w) \|_{{\cal H}^* \otimes {\cal H}^* \otimes {\mathbb R}^e} 
&=& 
2\hat{\mathbb E} [ 
\|\hat{D}^2 \Xi_2 (w, \,\cdot\,)_t \|_{{\cal H}^* \otimes {\cal H}^* \otimes {\mathbb R}^e}^2
]^{1/2}
\lesssim 
2\| \Xi_2 (w, \,\cdot\,)_t \|_{ {\mathbb D}_{2,2} }
\lesssim
\| \Xi_2 (w, \,\cdot\,)_t \|_{ L^{2} }.
\end{eqnarray*}

Next we take expectation in $w$-variable.  For any $r \ge 2$, 
\begin{eqnarray}
{\mathbb E}'  [\|D y_t  \|_{{\cal H}^* \otimes {\mathbb R}^e }^r]^{1/r}
&\lesssim&
{\mathbb E}' 
\Bigl[
\hat{\mathbb E} [   | \Xi_1 |^2]^{r/2}
\Bigr]^{1/r }
\lesssim
{\mathbb E}  [  | \Xi_1|^r]^{1/r}, 
\nn\\
%\qquad 
{\mathbb E}'  [ \|D^2 y_t  \|_{{\cal H}^* \otimes {\cal H}^*\otimes {\mathbb R}^e}^r]^{1/r}
&\lesssim&
{\mathbb E}' 
\Bigl[
\hat{\mathbb E} [   | \Xi_2 |^2]^{r/2}
\Bigr]^{1/r }
\lesssim
{\mathbb E}  [  | \Xi_2|^r]^{1/r}.
\nn
\end{eqnarray}
Therefore, if $\Xi_1, \Xi_2$ have moments of all order,
then we can show 
$y_t \in {\mathbb D}_{r,2} ({\mathbb  R}^e)$ for any $1 <r<\infty$.
(The Sobolev norm is non-decreasing both in $r$ and $k$.)

But, when one wants to make the above argument rigorous, 
the most difficult part is to define 
"something like stochastic integrals"
in (\ref{heu_xi1.eq}) and (\ref{heu_xi2.eq}).
In order to deal with this problem, 
we regard the system of ODEs (\ref{9heu_y.eq})--(\ref{heu_K.eq})
as 
the system of RDEs driven by the natural lift of $(w_t, b_t)$.
For instance, 
we understand RDE (\ref{9heu_y.eq}) in the following way;
\begin{equation}
\label{heu_y.eq}
%\nn
dy_t = \sum_{i=1}^d  V_i ( y_t) dw_t^i  + \sum_{i=1}^d  {\bf 0} \cdot db_t^i
\qquad
\qquad
\mbox{with \quad $y_0 =0\in {\mathbb R}^e$.}
\end{equation}
Here, ${\bf 0}$ is the constant vector fields that vanishes everywhere.
We understand
RDEs for $J$ and $K$ in the same way.

Then, we obtain 
$({\bf w}, {\bf b}; {\bf y}, {\bf J}, {\bf K})$
is a well-defined random rough path in ${\mathbb R}^{2d}
 \oplus {\mathbb R}^e \oplus {\rm Mat}(e,e)^{\oplus 2}$
and, by a similar argument,  
(any component of) any level path of this
has moments of all order.
If we take a projection to discard $b$-component, then 
we have the same $({\bf w}, {\bf y}, {\bf J}, {\bf K})$ as the one 
constructed from $w$ alone.

\begin{sloppypar}
The right hand side of (\ref{heu_xi1.eq}) can be interpreted 
as a rough path integral along $({\bf w}, {\bf b}; {\bf y}, {\bf J}, {\bf K})$
and consequently we get 
$({\bf w}, {\bf b}; {\bf y}, {\bf J}, {\bf K}, \bm{\Xi}_1)$.
Next, 
the right hand side of (\ref{heu_xi1.eq}) can be interpreted 
as a rough path integral along $({\bf w}, {\bf b}; {\bf y}, {\bf J}, {\bf K}, \bm{\Xi}_1)$.
So we get 
$({\bf w}, {\bf b}; {\bf y}, {\bf J}, {\bf K}, \bm{\Xi}_1, \bm{\Xi}_2)$.
The integrands (i.e., the ${\mathbb R}^e$-valued one-forms) in these integrations 
and  their derivatives are of at most polynomial growth.
Hence, 
(any component of) any level path of these rough paths
have moments of all order.
In particular, 
the first level paths of $\bm{\Xi}_1, \bm{\Xi}_2$ are $\Xi_1, \Xi_2$, respectively,
and hence they have moments of all order.
Thus, we have (formally) shown $y_t \in \cap_{1 <r <\infty} {\mathbb D}_{r,2} ({\mathbb R}^e)$.
\end{sloppypar}

To make the above argument rigorous, we consider piecewise linear approximations 
of $(w,b)$.
Then, all the integrals and ODEs are in the Riemann-Stieltjes sense
and it is not so difficult to prove ${\mathbb D}_{\infty}$-property 
for the approximating sequence.
After that we take limit and finish the proof.
Recall that
rough path theory is a very powerful tool to prove this kind of approximations,
thanks to Lyons' continuity theorem, etc.

%%%%%%%%%%%%%%%%%%%%%%%%%%%%%%%%%%%%%%%%%%%%%%%%%%%%%%%%%%%%%%%%%%%%%%
%\newpage
%
%
\section{Piecewise linear approximations}
\label{sec.pwl}
In this section we consider  ODEs driven by 
  piecewise linear approximations of $w$ and prove smoothness of the solution 
 in the sense of Malliavin calculus.
(We do not need advanced results in Malliavin calculus.
The first several chapters of \cite{nu} or \cite{sh} are enough.)
In this section
all  line integrals and ODEs are in the Riemann-Stieltjes sense.

Consider the following ODE 
driven by $x \in C_0^{1 -var} ([0,1], {\mathbb R}^d)$;
\begin{equation}
\label{rs_y.eq}
%\nn
dy_t = \sum_{i=1}^d  V_i ( y_t) dx_t^i 
\qquad
\qquad
\mbox{with \quad $y_0 =0\in {\mathbb R}^e$.}
\end{equation}
Write $y_t = I_t (x)$.
It is known that
It\^o map 
$I: C_0^{1 -var} ([0,1], {\mathbb R}^d) \to C^{1 -var} ([0,1], {\mathbb R}^e)$
is Fr\'echet smooth.
The Jacobian process and its inverse are given by 
\begin{eqnarray}
dJ_t &=& \sum_{i=1}^d  \nabla V_i ( y_t) \cdot J_t dx_t^i 
\qquad
\qquad
\mbox{with \quad $J_0 ={\rm Id}_e \in {\rm Mat}(e,e)$.}
\label{rs_J.eq}
\\
dK_t &=& - \sum_{i=1}^d K_t  \cdot \nabla V_i ( y_t)  dx_t^i 
\qquad
\qquad
\mbox{with \quad $K_0 ={\rm Id}_e \in {\rm Mat}(e,e)$.}
\label{rs_K.eq}
\end{eqnarray}
In the $1$-variation setting, we have a Gronwall-type  lemma.
So, it is not very difficult to show that
$ \|  J\|_{1 -var} + \|  K\|_{1 -var}\le C \exp( C \|  x\|_{1 -var})$
for some $C>0$ which is independent of $x$.

It is possible to write down the directional derivatives of $y_t$.
Let $h \in C_0^{1 -var} ([0,1], {\mathbb R}^d)$.
Below, we write $D_h^n y_t = (D_h)^n y_t = D^n I_t (x) \la h, \ldots, h\ra$ for brevity.
\begin{eqnarray}
 D_h y_t 
&=& 
J_t \int_0^t K_s  \sigma ( y_s) dh_s.
\label{rs_1int.eq}
\\
 D_h^2 y_t 
&=& 
J_t \int_0^t K_s \bigl\{ 
\nabla^2 \sigma ( y_s) \la D_h y_s, D_h y_s ,dw_s \ra
%\nn\\
%&&
%\qquad\qquad
+
2  \nabla \sigma ( y_s) \la D_h y_s ,dh_s\ra  
\bigr\}.
\label{rs_2int.eq}
\end{eqnarray}
For general $n =2,3,\ldots$,
\begin{eqnarray}
 D_h^n y_t 
&=& 
J_t \int_0^t K_s 
\Bigl\{
\sum_{l=2}^n   
\sum_{
i_1 + \ldots + i_l =n}
C_{i_1, \ldots, i_l}
\nabla^l \sigma ( y_s) \la D_h^{i_1} y_s ,  \ldots , D_h^{i_l} y_s, dw_s\ra  
\nn\\
&&
+
\sum_{l=1}^{n-1}   
\sum_{
i_1 + \ldots + i_l =n-1}
C'_{i_1, \ldots, i_l}
\nabla^l \sigma ( y_s) \la D_h^{i_1} y_s ,  \ldots , D_h^{i_l} y_s, dh_s\ra  
\Bigr\}
\label{rs_nint.eq}
\end{eqnarray}
Here, {\rm (i)}~ the summation $\sum_{i_1 + \ldots + i_l =n}$
runs over all non-decreasing sequence $0< i_1 \le  \ldots \le i_l$ of natural numbers
such that $i_1 + \ldots + i_l =n$,
{\rm (ii)}~
$C_{i_1, \ldots, i_l},  C'_{i_1, \ldots, i_l} \in {\mathbb N}$
are constants, 
but their exact values are not used in this paper.

The following remark is simple, but may be helpful.
\begin{re}\label{re.diag}
Let ${\cal X}$ and ${\cal Y}$ be real Banach spaces.
Suppose that $A$ and $B$ are two bounded, symmetric, $n$-multilinear maps 
from ${\cal X}^{\times n}$ to ${\cal Y}$.
If $A (x,x, \ldots,x) = B (x,x, \ldots,x)$ for any $x \in {\cal X}$,
then $A=B$ as multilinear maps.
(This can easily be checked as follows. 
Taking directional derivatives $D_{v_1} \cdots D_{v_n}$ of the map $x \mapsto A (x, \ldots,x)$,
we get $n! A (v_1, \ldots, v_n)$.)
\end{re}

%%%%%%%%%%%%\vspace{10mm}
%%%%%%%%%%%%

Let us consider the dyadic piecewise linear approximation of
the Gaussian process $w =(w_t)_{0 \le t \le 1}$. 
For $m =1,2, \ldots$, 
$w(m)=(w(m)_t)_{0 \le t \le 1}$ denotes the 
piecewise linear approximation associated with the partition $\{ l /2^m~|~ 0 \le l \le 2^m \}$ of $[0,1]$.
We write $\Delta^m_l w := w_{l /2^m} - w_{(l-1) /2^m}$ for simplicity.

Let $y(m)$ be the solution of (\ref{rs_y.eq}) with $x =w(m)$, i.e., $y(m) = I(w(m))$.
In a similar way, we define $J(m), K(m)$, etc.
These are functional on the abstract Wiener space $({\cal W}, {\cal H}, \mu)$.
It is intuitively clear that $y(m)_t$ is smooth in the sense of Malliavin calculus
 for each fixed $m$. 
However, we give a proof for completeness.

\begin{pr}\label{pr.app_smth}
For any $m=1,2,\ldots$ and $t \in [0,1]$, 
we have $y(m)_t \in {\mathbb D}_{\infty} ({\mathbb R}^e)$.
\end{pr}

\Proof
Since
$w(m)$ is a $C_0^{1 -var} ([0,1], {\mathbb R}^d)$-valued function 
of $\Delta^m_l w ~(1 \le l \le 2^m)$
and $I_t$ is Fr\'echet smooth,
there exists a smooth ${\mathbb R}^e$-valued function
$G =G_{t,m}$ defined on ${\mathbb R}^{d 2^m}$ 
such that
\[
y(m)_t = G (\Delta^m_1 w, \ldots,\Delta^m_{2^m} w).
\]
It is sufficient to show that 
$G$ and all of its partial derivatives are of at most exponential order.
(A vector-valued function $F$ defined 
on a Euclidean space is said to be of at most exponential order 
if there exists a constant $C>0$ such that
$|F(\xi)| \le C \exp (C |\xi|)$  for all $\xi$.)

First, note that $\|w(m)\|_{1 -var} \le \sum_{l} |\Delta^m_l w|$.
It is known that the solution $y$ of (\ref{rs_y.eq})
satisfies that $\|y\|_{1 -var} \le C (1 + \|x\|_{1 -var} )^C$ for some $C>0$.
Hence, $G$ itself 
is of at most polynomial growth in $\Delta^m_l w$'s.

%\begin{sloppypar}
%We can obtain 
% derivatives of $y(m)_t$ 
By replacing $w, h, y, J, K$ 
 with $w(m), h(m), y(m), J(m), K(m)$
in (\ref{rs_1int.eq})--(\ref{rs_nint.eq}), respectively,
we obtain explicit expressions of derivatives of $y(m)_t$. 
(Precisely, "polarization" is also needed.)
%\end{sloppypar}
%
The only quantities in those expressions that are not of  polynomial growth in $\|w(m)\|_{1 -var}$
is $\|J(m)\|_{1 -var}$ and $\|K(m)\|_{1 -var}$.
But, they are dominated by
$C \exp (C \|w(m)\|_{1 -var})$.
(Remark: Taking partial derivatives of $G$ corresponds to 
taking directional derivatives $D_{h}$ for some  $h$.
So, it suffices to estimate $D_{h_1} \cdots D_{h_n} y(m)_t$ for 
arbitrarily fixed $h_i$'s.)
%This completes the proof.
\QED

If $w$ is shifted by $h$, then $w(m)$ is shifted by $h(m)$. 
So, we have 
\begin{equation}
\label{m_D^n.eq}
D^n y(m)_t \la  h,h, \ldots, h\ra 
=
D^n I_t (w(m)) \la  h(m), h(m), \ldots, h(m)\ra, 
\end{equation}
where $D$ on the right side is in Fr\'echet sense, 
while 
$D$ on the left side is ${\cal H}$-derivative in Malliavin calculus.

%%%%%%%%%%%%\vspace{10mm}
%%%%%%%%%%%%

Let 
$(b_t) = (b_t^1, \ldots, b_t^d)$ be an independent copy of $(w_t)$.
The abstract Wiener space that corresponds to 
the $2d$-dimensional process $(w_t; b_t)$ is 
$({\cal W}^{\oplus 2}, {\cal H}^{\oplus 2}, \mu \times \mu)$.
The expectation with respect to $w$-variable and $b$-variable
are
denoted by ${\mathbb E}'$ and $\hat{\mathbb E}$, respectively.
The expectation with respect to $(w, b)$-variable 
is denoted by
${\mathbb E}= {\mathbb E}^{\mu \times \mu} = {\mathbb E}' \times \hat{\mathbb E}$.

%%%%%%%%%%%%\vspace{10mm}
%%%%%%%%%%%%

By replacing $h(m)$
on the right hand side of (\ref{m_D^n.eq}) with $b(m)$, we define
\begin{equation}\label{m.xi.def}
\Xi_n (m)_t (w,b) 
:=
D^n I_t (w(m)) \la  b(m), b(m), \ldots, b(m)\ra. 
\end{equation}
%
%For each $w$ and $b$, this is well-defined.
%
%
More explicitly, 
\begin{eqnarray}
 \Xi_1 (m)_t
&=& 
J(m)_t \int_0^t K(m)_s  \sigma ( y(m)_s) db(m)_s.
\label{m.xi_1int.eq}
%\nn
\end{eqnarray}
and, for $n=2,3,\ldots$,
\begin{eqnarray}
\lefteqn{
\Xi_n (m)_t
%}\nn\\
= 
J(m)_t \int_0^t K(m)_s
}
\nn\\
&&
\cdot
\Bigl\{
\sum_{l=2}^n   
\sum_{
i_1 + \ldots + i_l =n}
C_{i_1, \ldots, i_l}
\nabla^l \sigma ( y(m)_s) \la  \Xi_{i_1} (m)_s ,  \ldots ,  \Xi_{i_l} (m)_s, dw(m)_s\ra  
\nn\\
&&
\qquad+
\sum_{l=1}^{n-1}   
\sum_{
i_1 + \ldots + i_l =n-1}
C'_{i_1, \ldots, i_l}
\nabla^l \sigma ( y(m)_s) \la  \Xi_{i_1} (m)_s,  \ldots , \Xi_{i_l} (m)_s, db (m)_s\ra  
\Bigr\},
\label{m.xi_n.eq}
\end{eqnarray}
where the constants are the same as in (\ref{rs_nint.eq}).

%%%%%%%%%%%%\vspace{10mm}
%%%%%%%%%%%%

Since $\Xi_n (m)_t (w,b)$  is defined for all $w$ and $b$,
we can think of $\Xi_n (m)_t (w,\,\cdot\,)$ as a Wiener functional in $b$ for each fixed $w$.
Then, it is clear from the right hand side of (\ref{m.xi.def})
that, for each $w$,  $\Xi_n (m)_t (w,b)$ is a polynomial of order $n$
in $\Delta^m_l b$'s.
In paticular,  
$\Xi_n (m)_t (w,\,\cdot\,)$ belongs to $n$th order inhomogeneous Wiener chaos.
Moreover, by straight-forward computation,
\[
\hat{D}^n \Xi_n (m)_t (w,\,\cdot\,) \la h, \ldots, h \ra
=
n! D^n y(m)_t \la  h,h, \ldots, h\ra 
\qquad
(h \in {\cal H}).
\]
Here, $\hat{D}$ stands for ${\cal H}$-derivative with respect to $b$-variable.
Note that both side do not depend on $b$.
It follows from Remark \ref{re.diag} that 
$\hat{D}^n \Xi_n (m)_t (w,\,\cdot\,) = n! D^n y(m)_t \in ({\cal H}^*)^{\otimes n}$
for each $w$.

The next proposition implies that, if $\{ \Xi_n(m)_t \}$ is Cauchy in $L^r$-norm, 
then $\{D^n y(m)_t\}$ is also Cauchy  in $L^r$-norm.
As a result, the proof of the main theorem reduces to showing 
 $\{ \Xi_n(m)_t \}$ is Cauchy in $L^r$-norm for any $t, n,$ and  $r \in [2, \infty)$.
\begin{pr}\label{pr.est_ym}
For $n=1,2,\ldots$ and $2 \le r<\infty$,
there is a positive constant $C=C_{r,n}$ such that
\[
{\mathbb E}'  [ \|D^n y(m)_t  \|_{{\cal H}^{* \otimes n} \otimes {\mathbb R}^e}^r]^{1/r}
\le 
C  {\mathbb E}  [  | \Xi_n(m)_t |^r]^{1/r}
\]
for any $m=1,2,\ldots$ and $0 \le t \le 1$. 
In a similar way, we have 
\[
{\mathbb E}'  [ \|D^n y(m)_t - D^n y(m')_t  \|_{{\cal H}^{* \otimes n}\otimes {\mathbb R}^e}^r]^{1/r}
\le 
C  {\mathbb E}  [  | \Xi_n(m)_t - \Xi_n(m' )_t |^r]^{1/r}
\]
for any $m, m' =1,2,\ldots$ and $0 \le t \le 1$. 
\end{pr}

\Proof
We prove the first assertion.
Note that 
$$
\|D^n y(m)_t  \|_{{\cal H}^{* \otimes n} \otimes {\mathbb R}^e} 
= 
\frac{1}{n!} \| \hat{D}^n \Xi_n (m)_t (w,\,\cdot\,) \|_{{\cal H}^{* \otimes n} \otimes {\mathbb R}^e} 
 =
 \frac{1}{n!}  
 \hat{\mathbb E} [ 
  \| \hat{D}^n \Xi_n (m)_t (w,\,\cdot\,) \|_{{\cal H}^{* \otimes n} \otimes {\mathbb R}^e}^2
  ]^{1/2}.
  $$
Here, we used the fact that
$\hat{D}^n \Xi_n (m)_t (w,\,\cdot\,)$ does not depend on $b$.
Since all ${\mathbb D}_{2, n}$-norms are equivalent on a fixed inhomogeoneous
Wiener chaos ($n=0,1,2,\ldots$),  
the right hand side is dominated by 
$C_{r,n} \hat{\mathbb E} [  | \Xi_n (m)_t (w,\,\cdot\,) |^2]^{1/2}$.
Taking $L^r$-norm of this inequality with respect to $w$-variable, 
we show the first assertion. We can show the second one in a similar way.
\QED

%%%%%%%%%%%%%%%%%%%%%%%%%%%
%\newpage
%%%%%%%%%%%%%%%%%%%%%%%%%%

\section{Rough path theory}
In this section, we recall basic results in rough path theory.
We use T. Lyons' original formulation of rough path integral and rough differential equation
as in Lyons-Qian \cite{lq}, Lyons-Caruana-L\'evy \cite{lcl}. (See also Lejay \cite{le}.)
Throughout this paper  we assume $2 \le p <4$, where $p$ denotes the roughness constant.

%%%%%%%%%%%%%%%%%
\subsection{Deterministic operations on rough path spaces}

In this subsection
we summarize various deterministic operations in rough path theory,
which will be needed in what follows.
No probability measures or random variables appear in this subsection.

Let $G\Omega_p ({\mathbb R}^d)$ be the geometric rough path space 
with $p$-variation topology.
A generic element of $G\Omega_p ({\mathbb R}^d)$
is denoted by ${\bf x } = ({\bf x }^1, \ldots, {\bf x }^{[p]})$.
It satisfies an algebraic relation called Chen's identity 
and $i$th level path ${\bf x }^{i}$ is of finite $p/i$-variation for $1 \le i \le [p]$.
The intrinsic control function of ${\bf x } \in G\Omega_p ({\mathbb R}^d)$ is denoted by
$\omega_{{\bf x }} (s,t) := \sum_{i=1}^{[p]} \| {\bf x }^{i} \|_{p/i -var, [s,t]}^{p/i}$,
where 
$\| \,\cdot\, \|_{p -var, [s,t]}$ stands for $p$-variation norm 
restricted on the subinterval $[s,t] \subset [0,1]$.

If $f : {\mathbb R}^d \to {\rm Mat}(e,d)$ is $C^{[p]+1}$, 
then the rough path integral 
$\int f( {\bf x }) d {\bf x } \in G\Omega_p ({\mathbb R}^e)$ is well-defined
and extends Riemann-Stieltjes integral.
Moreover, this integration map is continuous with respect to the rough path topology.
(See \cite{lq} or \cite{lcl} for these facts.  We call $f$  the integrand.)
We often take $\hat{f} = {\rm Id}_d \oplus f$ as the integrand.
Then, $({\bf x}, \int f( {\bf x }) d {\bf x })$
is a well-defined element in $ G\Omega_p ({\mathbb R}^d \oplus {\mathbb R}^e )$.
We say the integrand $f$ is at most of polynomial growth if 
there exists $c>0$ such that 
$|\nabla^j f (\xi) | \le c (1 +|\xi| )^c$ holds for all $\xi \in {\mathbb R}^d$
and $0 \le j \le [p]+1$.
(Though it is probably better to say "$f$ is at most of polynomial growth with its derivatives,"
we use this terminology for simplicity.)
If $f$ at most of polynomial growth, then so is $\hat{f}$.

%%%%%%%%%%%%%\vspace{10mm}
%%%%%%%%%%%%%%

Now we recall a growth estimate of solutions of RDEs in a general setting.
Let $\sigma : {\mathbb R}^e \to {\rm Mat}(e,d)$
be of $C_b^{[p] +1}$ (i.e., $M : = \sum_{j=0}^{[p]+1} \| \nabla^j \sigma \|_{\infty}<\infty$).
We consider RDE
$dy_t = \sigma (y_t) dx_t$ with a given $y_0$.
When ${\bf x}$ is the driving rough path, then  a unique solution
${\bf z} =({\bf x}, {\bf y})$ satisfies the following estimate:
\begin{eqnarray*}
 |{\bf z}^i_{s,t}|
 &\le&
c (1+M)^c  (1 + \omega_{{\bf x}}(0,1))^c   \omega_{{\bf x}} (s,t)^{i/p} 
\qquad
(0 \le s \le t \le 1, ~1 \le i \le [p]),
\end{eqnarray*}
where $c>0$ is a constant independent of ${\bf x}, M$, which may change from line to line.
Hence, we have 
\begin{equation}
\| {\bf z}^i \|_{p /i -var}
\le 
c (1+M)^c 
 (1 + \sum_{i=1}^{[p]} \| {\bf x}^i \|_{p/i -var}^{p/i} )^{c}
\qquad
(1 \le i \le [p]).
\label{z_hyoka.ineq}
\end{equation}

Next we consider the case where $\sigma$ is of $C^{[p] +1}$,
but $\sigma$ and its derivatives may have linear growth, that is, 
$\sum_{j=0}^{[p]+1} | \nabla^j \sigma (\xi) | \le c (1 +|\xi|)$ for all $\xi \in {\mathbb R}^e$.
In this case, the RDE may not have a global solution.
(If it exists, then it is unique.)
So, assume that a global solution ${\bf z}$ exists for ${\bf x}$.
We set $y_t = y_0 +{\bf y}^1_{0,t}$.
Since $y$ stays inside the ball of radius $\|y\|_{\infty}$ centered at the origin,
we only use information of $\sigma$ restricted on the ball of radius $2\|y\|_{\infty}$.
So, we may take $M =  2c (1 + \|y\|_{\infty})$ in (\ref{z_hyoka.ineq}).
(Use a cutoff argument).
Thus, we obtain
\begin{equation}
\| {\bf z}^i \|_{p /i -var}
\le 
c (1+  \|y\|_{\infty})^c 
 (1 + \sum_{i=1}^{[p]} \| {\bf x}^i \|_{p/i -var}^{p/i} )^{c}
\qquad
(1 \le i \le [p]).
\label{z_hyoka2.ineq}
\end{equation}
for this linear growth case if a global solution exists.

%%%%%%%%%%%%%\vspace{10mm}
%%%%%%%%%%%%%%

Now we give a remark on geometric rough paths on 
the direct sum of two (or more) vector spaces. 
An element in $G\Omega_p ({\mathbb R}^d \oplus {\mathbb R}^e )$
is sometimes written as 
$({\bf x}, {\bf y})$, where each  "component" is an element of  
$G\Omega_p ({\mathbb R}^d)$ or $G\Omega_p ({\mathbb R}^e)$, respectively.
This is slightly misleading because it looks like an element of 
$G\Omega_p ({\mathbb R}^d) \times G\Omega_p ({\mathbb R}^e)$
and "cross terms" of high level paths may be forgotten.
However, we use this kind of notation when risk of confusion is low,
simply because we do not know a better way.
Typical examples are as follows;
\\
\\
\noindent
{\rm (i)}~In Lyons' original formulation, a solution ${\bf z}$ of RDE 
$y = \sum_i V_i (y) dx^i$ driven by ${\bf x} \in G\Omega_p ({\mathbb R}^d) $
 is actually an element in 
$G\Omega_p ({\mathbb R}^d \oplus {\mathbb R}^e )$.
But, we often write ${\bf z} = ({\bf x}, {\bf y})$,
where ${\bf y} \in G\Omega_p ({\mathbb R}^e)$ is the image of ${\bf z}$ 
by the projection map onto the second component.
(A projection map in the sense of linear algebra naturally extends 
to a projection map  in the sense of rough path.)
The same remark also goes for $({\bf x}, {\bf y}, {\bf J}, {\bf K})$, etc.
\\
\noindent
{\rm (ii)}~As is explained above,  
$({\bf x}, \int f( {\bf x }) d {\bf x })$ should be understood as 
rough path integral along ${\bf x}$ against the integrand $\hat{f} = {\rm Id}_d \oplus f$.
\\
\noindent
{\rm (iii)}~In the previous sections we denote the $2d$-dimensional Gaussian process 
by $(w_t, b_t)$.
Corresponding to this notation, we denote a generic element of  
$G\Omega_p ({\mathbb R}^{d} \oplus {\mathbb R}^{d}  )$ by $({\bf x}, {\bf v})$.
(The reader may find this one
 a bit unnatural, but we believe this notation will turn out to be useful later.)
\\
\\

We give another remark for the linear combination 
and the multiplication of two components of a rough path.
For example,  for $({\bf x}, {\bf x}' ) \in G\Omega_p ({\mathbb R}^d \oplus {\mathbb R}^d )$,
the linear combination  $a{\bf x} +b {\bf x}' \in G\Omega_p ({\mathbb R}^d)$
is well-defined since 
it can be understood as 
 rough path integral $\int \{a d{\bf x} + b d{\bf x}' \}$.
Similarly, 
for $({\bf x}, {\bf M} ) \in G\Omega_p ({\mathbb R}^d \oplus {\rm  Mat}(e,d) )$,
${\bf M} \cdot {\bf x} \in G\Omega_p ({\mathbb R}^e)$
can be understood in the same way since 
${\bf M} \cdot {\bf x}  = \int\{  {\bf M} \cdot d {\bf x} + d {\bf M} \cdot  {\bf x}\}$.
Here, the "dot" stands for the matrix multiplication.
Therefore, these are special cases of rough path integration
against an integrand of at most polynomial growth.
To keep our exposition concise, 
we will not treat these cases independently from now on.

%%%%%%%%%%%%%\vspace{10mm}
%%%%%%%%%%%%%%

For a given $({\bf x}, {\bf v}) \in G\Omega_p ({\mathbb R}^{2d})$,
we consider a system of RDEs (\ref{rs_y.eq})--(\ref{rs_K.eq}). 
(The coefficients of $dv^i$ ($1 \le i \le d$) are simply ${\bf 0}$.
See (\ref{heu_y.eq}).)
Consequently, 
\[
G\Omega_p ({\mathbb R}^{2d})
\ni  ({\bf x}, {\bf v})  \mapsto   ({\bf x}, {\bf v}; {\bf y}, {\bf J}, {\bf K})
\in 
G\Omega_p ({\mathbb R}^{2d}  \oplus {\mathbb R}^e \oplus {\rm Mat} (e,e)^{ \oplus 2} )
\]
is locally Lipschitz continous with respect to $p$-variation distance.
Note that
"$(x,y,J,K)$-component" of the above rough path actually depends only on ${\bf x}$.
(In other words, 
if we dicard 
"$v$-component"  from $ ({\bf x}, {\bf v}; {\bf y}, {\bf J}, {\bf K})$ 
by using a suitable projection, 
then we get the same 
$ ({\bf x},  {\bf y}, {\bf J}, {\bf K})$ which is constructed from ${\bf x}$ alone.)
When $({\bf x}, {\bf v})$ is the natural lift of 
$(x, v) \in C_0^{1 -var} ([0,1], {\mathbb R}^{2d})$,
then the first level of 
$ ({\bf x}, {\bf v}; {\bf y}, {\bf J}, {\bf K})$ coincides with the solution 
$(x,v;y,J,K)$ of ODEs in Riemann-Stieltjes sense (after the initial value is adjusted,
which means 
 $J_t = {\rm Id} + {\bf J}^1_{0,t}$ and $K_t = {\rm Id} + {\bf K}^1_{0,t}$).

%%%%%%%%%%%%%\vspace{10mm}
%%%%%%%%%%%%%%

Next we define the rough path extension of $\Xi_n$ 
so that the relation
$\bm{\Xi}_n ({\bf x},{\bf v})= (D^n I_t )({\bf x}) \la {\bf v}, \ldots, {\bf v} \ra$
still
holds formally.
Below, vector spaces ${\cal X}_n ~(n=0,1,2,\ldots)$ are defined as follows:
$
{\cal X}_0:={\mathbb R}^{2d}  \oplus {\mathbb R}^e \oplus {\rm Mat} (e,e)^{ \oplus 2}
$
%\qquad
%\mbox{and}
%\qquad
and
$
{\cal X}_n := {\cal X}_0 \oplus ({\mathbb R}^e)^{\oplus n}
$
%\qquad
%\mbox{
for $n \ge 1.$
%$}\]

Let us denote an element of $G\Omega_p ( {\cal X}_0)$
by
$ ({\bf x}, {\bf v}; {\bf y}, {\bf J}, {\bf K})$.
We define $\bm{\Xi}_1$ as follows:
\begin{equation}
\bm{\Xi}_1
=
({\rm Id} + {\bf J}) 
\cdot \int  ({\rm Id} + {\bf K})  \sigma ({\bf y}) d{\bf v}.
\label{rp_xi1.eq}
\end{equation}
This is a rough path integral along $({\bf x}, {\bf v}; {\bf y}, {\bf J}, {\bf K})$
against a polynomially growing integrand.
Hence, the map 
$$
G\Omega_p ( {\cal X}_0)  \ni ({\bf x}, {\bf v}; {\bf y}, {\bf J}, {\bf K}) 
 \mapsto ({\bf x}, {\bf v}; {\bf y}, {\bf J}, {\bf K}, \bm{\Xi}_1)  \in G\Omega_p ( {\cal X}_1)
  $$
is well-defind and continuous.

We continue this procedure recursively.
For $ ({\bf x}, {\bf v}; {\bf y}, {\bf J}, {\bf K}, \bm{\Xi}_1, \ldots, \bm{\Xi}_{n-1}) 
\in G\Omega_p ( {\cal X}_{n-1})$,
\begin{eqnarray}
\bm{\Xi}_n
&=& 
({\rm Id} + {\bf J}) 
\cdot \int ({\rm Id} + {\bf K}) 
\Bigl\{
\sum_{l=2}^n   
\sum_{
i_1 + \ldots + i_l =n}
C_{i_1, \ldots, i_l}
\nabla^l \sigma ( {\bf y}) \la \bm\Xi_{i_1} ,  \ldots , \bm\Xi_{i_l}, d{\bf w}\ra  
\nn\\
&&
+
\sum_{l=1}^{n-1}   
\sum_{
i_1 + \ldots + i_l =n-1}
C'_{i_1, \ldots, i_l}
\nabla^l \sigma ({\bf y}) \la  \bm\Xi_{i_1},  \ldots , \bm\Xi_{i_l}, d{\bf v}  \ra  
\Bigr\}.
\label{rp_xin.eq}
\end{eqnarray}
Note that (\ref{rp_xi1.eq}) and (\ref{rp_xin.eq})
are parallel to (\ref{rs_1int.eq}) and (\ref{rs_nint.eq}), respectively.
Again this is a rough path integral along 
against a polynomially growing integrand.
Hence, the map 
$$
G\Omega_p ( {\cal X}_{n-1})  \ni ({\bf x}, {\bf v}; {\bf y}, {\bf J}, {\bf K}, \bm{\Xi}_1, \ldots \bm{\Xi}_{n-1}) 
 \mapsto ({\bf x}, {\bf v}; {\bf y}, {\bf J}, {\bf K}, \bm{\Xi}_1, \ldots \bm{\Xi}_n)  \in G\Omega_p ( {\cal X}_n)
  $$
is well-defind and continuous for all $n \ge 2$.

%%%%%%%%%%%%%\vspace{10mm}
%%%%%%%%%%%%%%

In summary, we have the following sequence of continuous maps
between geometric rough path spaces:
\begin{eqnarray}
({\bf x}, {\bf v})  \in  G\Omega_p ({\mathbb R}^{d}  \oplus {\mathbb R}^{d}  )
&\mapsto&
({\bf x}, {\bf v}; {\bf y}, {\bf J}, {\bf K})
\in 
G\Omega_p ({\cal X}_0)
\nn\\
&\mapsto&
({\bf x}, {\bf v}; {\bf y}, {\bf J}, {\bf K}, \bm{\Xi}_1)
\in 
G\Omega_p ({\cal X}_1)
\nn\\
&\mapsto& \cdots\cdots
\nn\\
&\mapsto&
 ({\bf x}, {\bf v}; {\bf y}, {\bf J}, {\bf K}, \bm{\Xi}_1, \ldots \bm{\Xi}_{n-1})  
 \in G\Omega_p ( {\cal X}_{n-1})
 \nn\\
 &\mapsto&
 ({\bf x}, {\bf v}; {\bf y}, {\bf J}, {\bf K}, \bm{\Xi}_1, \ldots \bm{\Xi}_n)  \in G\Omega_p ( {\cal X}_n)
\nn\\
&\mapsto& \cdots\cdots.
\label{seq_map}
 \end{eqnarray}

We end this subsection with the following remark.
\begin{re}\label{re.mxi}
Let $({\bf w}(m), {\bf b}(m))$ be the natural lift of 
the dyadic piecewise linear approximation $(w(m), b(m))$.
If we take $({\bf w}(m), {\bf b}(m))$ as the input $({\bf x}, {\bf v})$ in (\ref{seq_map}),
then the first level path of 
$\bm{\Xi}_n$ coincides with $\Xi_n (m)(w,b)$ 
in the previous section for all $n$.
\end{re}

%%%%%%%%%%%%%\vspace{10mm}
%%%%%%%%%%%%%%\newpage
\subsection{Some probabilistic results on rough path space}

In this subsection we present some basic
probabilistic results on rough path space.
For a while, random variables are defined on an arbitrary probability space $(\Omega, {\mathbb P})$.
We often write $L^{\infty -} :=\cap_{1 <r< \infty} L^r$.
Recall that the $i$th level path of a geometric rough path is an element of a Banach space
$C^{p/i -var} (\triangle, ({\mathbb  R}^{d})^{\otimes i})$,
the space of continous maps 
from $\triangle =\{(s,t)~|~0 \le s \le t \le 1 \}$ to $({\mathbb  R}^{d})^{\otimes i}$
with finite $p/i$-variation.

\begin{df}
{\rm (i)}~Let $Z_m~(m=1,2,\ldots)$ be $L^{\infty -}$-random variables 
that takes values in a real Banach space ${\cal B}$.
We say $\{Z_m\}$ is bounded in $L^{\infty -}$ if it is bounded in $L^r$
for any $r \in (1, \infty)$.
\\
{\rm (ii)}~
let  $Z_m~(m=1,2,\ldots, \infty)$ be as above.
We say $\{Z_m\}$ converges to $Z_{\infty}$ in $L^{\infty -}$ 
if it converges in $L^r$ for any $r \in (1, \infty)$.
\\
{\rm (iii)}~
Each level path of a $G\Omega_p ({\mathbb R}^{d})$-valued random variable
is a Banach space-valued random variable.
So we use the same terminologies for (a sequence of)
$G\Omega_p ({\mathbb R}^{d})$-valued random variables.
\end{df}

The following lemma will turn out to quite useful.
Thanks to this lemma, we need not estimate $L^r$-norm of difference between
two rough path space-valued 
random variables.
\begin{lm}
\label{lm.uni_integ}
Let $Z_m~(m=1,2,\ldots, \infty)$ be $L^{\infty -}$-random variables 
that takes values in a real Banach space ${\cal B}$,
which are bounded in $L^{\infty -}$.
Assume further that $Z_m$ converges to $Z_{\infty}$ a.s.
Then, $\lim_{m \to \infty} Z_m =Z_{\infty}$ in $L^{\infty -}$. 
\end{lm}

\Proof
For $r \in (1,\infty)$,
we set $Y_m = \|Z_m -Z\|_{{\cal B}}^r$, which is  real-valued.
Obvoiusly, 
the $L^2$-norm of $Y_m$ is bounded. 
So, $\{Y_m\}$ is uniformly integrable.
Since $Y_m \to 0$ a.s., $Y_m$ converges to $0$  in $L^1$, which means 
$\lim_{m \to \infty} Z_m =Z_{\infty}$ in $L^{r}$ for any $r$.
\QED

%%%%%%%%%%%%%\vspace{10mm}
%%%%%%%%%%%%%%

\begin{lm}\label{lm.conv_rpi}
Let ${\bf z}_m ~(m=1,2,\ldots, \infty)$ be $G\Omega_p ({\mathbb R}^d)$ valued-random variables.
We assume  
$\{ {\bf z}_m\}$ 
is bounded in $L^{\infty -}$ and it converges to ${\bf z}_{\infty}$ a.s.
For a $C^{[p]+1}$-one form 
$f : {\mathbb R}^d \to {\rm Mat}(e,d)$, which is of at most polynomial growth,
we set 
${\bf a}_m = \int f ( {\bf z}_m) d{\bf z}_m$.
Then, the sequence 
$\{ {\bf a}_m\}$ of $G\Omega_p ({\mathbb R}^d)$ valued-random variables
converges to ${\bf a}_{\infty}$ in $L^{\infty -}$ and a.s.
\end{lm}

\Proof 
First, suppose that $M := \sum_{j=0}^{[p]+1} \| \nabla^j f \|_{\infty}<\infty$.
If ${\bf x}$ is a geometric rough path controlled by a control function $\omega$, i.e., 
\[
|{\bf x}^i_{s,t}| \le \omega(s,t)^{i/p} 
\qquad
\qquad
(0 \le s \le t \le 1, \quad 1 \le i \le [p]),
\]
then there is $c>0$ which is independent of ${\bf x}, \omega, M, (s,t)$
such that
\[
 |{\bf a}^i_{s,t}|
 \le
c M^i  (1 + \omega(0,1))^c \omega(s,t)^{i/p} 
\qquad
(0 \le s \le t \le 1, \quad 1 \le i \le [p]),
\]
where ${\bf a} = \int f ( {\bf x}) d{\bf x}$.
We can choose the intrinsic control $\omega_{{\bf x}}$ 
as the control function. 
Hence, we have 
$$
\| {\bf a}^i \|_{p /i -var}
 \le c M^i  (1 + \omega_{{\bf x}}(0,1))^{c +i/p}
\le 
c M^i  (1 + \sum_{i=1}^{[p]} \| {\bf x}^i \|_{p/i -var}^{p/i} )^{c +i/p}.
$$

Next we consider the polynomially growing case.
The first level path $x_t = {\bf x}^1_{0,t}$
stays inside a ball of radius $\|x\|_{\infty}$.
Hence, we only use information of $f$ 
restricted on the ball of radius $2\|x\|_{\infty}$.
(Use the cutoff argument if necessary.)
Therefore, 
we may take $M = c_1 ( 1 + \|x\|_{\infty})^{c_1}$ and we have
\begin{eqnarray}
\| {\bf a}^i \|_{p /i -var}
 &\le& c ( c_1 ( 1 + \|x\|_{\infty})^{c_1} )^i  
(1 + \sum_{i=1}^{[p]} \| {\bf x}^i \|_{p/i -var}^{p/i} )^{c +i/p} 
 \nn\\
 &\le&
 c_2   
(1 + \sum_{i=1}^{[p]} \| {\bf x}^i \|_{p/i -var} )^{c_2} 
\label{f_int.ineq}
\end{eqnarray}
for some $c_2>0$ which is independent of ${\bf x}$.  
%
%Note that this is a deterministic estimate.

It follows
from the deterministic estimate (\ref{f_int.ineq})
that
$\{ {\bf a}_m\}$ is bounded in $L^{\infty -}$.
Since rough path integration map is continuous, 
it is clear that ${\bf a}_m \to {\bf a}_{\infty}$ a.s. as $m\to \infty$. 
Using Lemma \ref{lm.uni_integ}, we finish the proof.
\QED

%%%%%%%%%%%%%\vspace{10mm}
%%%%%%%%%%%%%%

From now on, the probability space is 
$({\cal W} \oplus {\cal W},  \mu \times \mu)$.
A generic element of ${\cal W} \oplus {\cal W}$ is denoted by $(w, b)$ as before.
$(w(m), b(m))$ stands for the $m$th dyadic piecewise linear approximation of $(w, b)$ and
its natural lift is denoted by $({\bf w}(m), {\bf b}(m))$.

Assume  that the covariance $R(s,t) = {\mathbb E}' [w^1_s w^1_t]$ 
is of finite 2D $\rho$-variation for some $\rho \in [1,2)$.
Then,
$\{ ({\bf w}(m), {\bf b}(m)) \}_{m=1}^{\infty}$
converges 
to some $G\Omega_p ({\mathbb R}^{2d})$-valued random variable $({\bf w}, {\bf b})$
in $L^{\infty -}$ if $2 \rho < p <4$.
(see Theorem 15.42, p. 436, \cite{fvbk}).
%
%Moreover, ${\cal H}$ is embedded in $C_0^{\rho -var} ([0,1], {\mathbb R}^d)$ under this assumption.
%
When $w$ is fractional Brownian motion with $H \in (1/4, 1/2]$,
then  $2\rho =1/H$. (See Proposition 15.5, \cite{fvbk}.)

By taking a subsequence if necessary, 
we may assume that $\lim_{m \to \infty} ({\bf w}(m), {\bf b}(m)) = ({\bf w}, {\bf b})$ a.s.
Abusing notation, we denote the subsequence by $\{ ({\bf w}(m), {\bf b}(m)) \}$ again.
(If the covariance $R(s,t)$ is "H\"older dominated," then taking a subsequence is not necessary.
See Excercise 15.44,  \cite{fvbk}.)
$({\bf w}, {\bf b})$ is called the natural lift of $(w, b)$.

\begin{df}\label{def_wbm}
Consider the sequence of continuous maps in (\ref{seq_map}).
For the input $({\bf x}, {\bf v}) =({\bf w}, {\bf b})$,
the output is simply denoted by 
$({\bf w}, {\bf b}; {\bf y}, {\bf J}, {\bf K}, \ldots \bm{\Xi}_n)$.
For the input $({\bf x}, {\bf v}) =({\bf w}(m), {\bf b}(m))$,
the output is denoted by 
$({\bf w}(m), {\bf b}(m); {\bf y}(m), {\bf J}(m), {\bf K}(m), \ldots \bm{\Xi}_n(m))$.
We supress the dependency on $(w,b)$ or $({\bf w}, {\bf b})$ to keep the notations simple.
\end{df}

%%%%
The following is a key technical lemma.
It is a slight modification of Cass-Litterer-Lyons' integrability lemma in \cite{cll}
on integrability of Jacobian process.
The point is uniformity in $m$. 
(This lemma holds for the original sequence, 
not just for the a.s. convergent subsequence we have chosen.)
\begin{lm}\label{lm.unif_exp_int}
Assume {\bf (H)}.~
We use the same notation as in Definition \ref{def_wbm} above.
We set $J(m)_t ={\rm Id} + {\bf J(m)}^1_{0,t}$, $K(m)_t ={\rm Id} + {\bf K}(m)^1_{0,t}$, etc.
Then, for any $r \in (1, \infty)$, we have
\[
\sup_m {\mathbb E}[ \|J(m)\|_{\infty}^r ] + \sup_m  {\mathbb E}[ \|K (m)\|_{\infty}^r ]
+
{\mathbb E}[ \|J\|_{\infty}^r ] + {\mathbb E}[ \|K\|_{\infty}^r ] <\infty.
\]
\end{lm}

\Proof
The proof will be given in Section \ref{sec.pr.lm} below.
\QED

%%%%%%%%%%%%%\vspace{10mm}
%%%%%%%%%%%%%%

Now we work with the a.s. converging subsequence again.

\begin{lm}\label{lm.conv_wtoK}
We use the same notation as in Definition \ref{def_wbm} above.
Then,
\[
({\bf w}(m), {\bf b}(m); {\bf y}(m), {\bf J}(m), {\bf K}(m) ) 
\to 
({\bf w}, {\bf b}; {\bf y}, {\bf J}, {\bf K})
\qquad
\mbox{as $m \to \infty$}
\] 
in $L^{\infty -}$ and a.s.
\end{lm}

\Proof
Consider the system of RDEs (\ref{rs_y.eq})--(\ref{rs_K.eq}) driven by $({\bf x}, {\bf v})$.
Although the coeficient is not  $C_b^{[p]+1}$,
it is known that 
 (\ref{rs_y.eq})--(\ref{rs_K.eq})  have a global solution for 
 any $({\bf x}, {\bf v}) \in G\Omega_p ({\mathbb R}^{2d})$
and, moreover, Lyons' continuity theorem holds.
Then, almost sure convergence follows immediately.
By Lemma \ref{lm.uni_integ}, it suffices to show 
$\{ ({\bf w}(m), {\bf b}(m); {\bf y}(m), {\bf J}(m), {\bf K}(m) )  \}_{m=1}^{\infty}$ 
is bounded in $L^{\infty -}$.

For a while, our argument will be deterministic.
Let us denote by $\hat \sigma$ the coefficient of 
RDEs (\ref{rs_y.eq})--(\ref{rs_K.eq}), which is a $C^{[p]+1}$-function on 
${\mathbb R}^e \oplus {\rm Mat} (e,e)^{\oplus 2}$.
By  abusing notations, we write a generic element of this set by $(y, J, K)$.  
Then, 
$\hat \sigma$ and its derivatives are linear growth in $J$ and $K$,
but bounded in $y$.
Precisely, there exists $c>0$ such that
\begin{equation}
\sum_{j=0}^{[p]+1} | \nabla^j \hat\sigma (y,J,K) | \le c (1 +|J| +|K|)
\qquad
\mbox{for all $(y,J,K) \in {\mathbb R}^e \oplus {\rm Mat} (e,e)^{\oplus 2}.$ }
\nn
\end{equation}
Hence, we can use a slight modification of (\ref{z_hyoka2.ineq})
to obtain the following deterministic estimate for $1 \le i \le [p]$;
\begin{equation}
\|  ( {\bf x}, {\bf v}; {\bf y}, {\bf J}, {\bf K}  )^i \|_{p /i -var}
\le 
c (1+  \|J\|_{\infty} +  \|K\|_{\infty} )^c 
 (1 + \sum_{i=1}^{[p]} \| ( {\bf x}, {\bf v} )^i \|_{p/i -var}^{p/i} )^{c}
%
%\quad
%(1 \le i \le [p]).
%
\label{z_hyoka3.ineq}
\end{equation}
Here, $c>0$ is independent of $( {\bf x}, {\bf v} ) \in G\Omega_p ({\mathbb R}^{2d})$.

Set $( {\bf x}, {\bf v} ) = ( {\bf w} (m), {\bf b} (m))$ in (\ref{z_hyoka3.ineq}).
Then, the second factor on the right hand side is clearly  bounded in $L^{\infty -}$. 
So is the first factor by Lemma \ref{lm.unif_exp_int}.
This completes the proof of the lemma.
\QED

%%%%%%%%%%%%%
%\vspace{10mm}
%%%%%%%%%%%%%%

By just combining the results we have already proved, 
we obtain the following lemma:
\begin{lm}\label{pr.conv_all}
We use the same notation as in Definition \ref{def_wbm} above.
Then, for any $n \ge 1$,
\begin{eqnarray*}
\lefteqn{
({\bf w}(m), {\bf b}(m); {\bf y}(m), {\bf J}(m), {\bf K}(m), \bm{\Xi}_1(m) ,\ldots, \bm{\Xi}_n(m)) 
}
\\
&\to& 
({\bf w}, {\bf b}; {\bf y}, {\bf J}, {\bf K}, \bm{\Xi}_1 ,\ldots, \bm{\Xi}_n)
\qquad
\mbox{as $m \to \infty$ in $L^{\infty -}$ and a.s.}
\end{eqnarray*}
%in $L^{\infty -}$ and a.s.
\end{lm}

\Proof
It immediately follows from (\ref{seq_map}),
Lemma \ref{lm.conv_rpi}, and Lemma \ref{lm.conv_wtoK}
\QED

%%%%%%%%%%%%%\vspace{10mm}
%%%%%%%%%%%%%%

\subsection{Proof of main theorem}\label{subsec.proof}

Now we are in a position to prove our main theorem in this paper (Theorem \ref{tm.main}).
First, note that ${\bf y}^1$ depends only on $w$-componet of $(w, b)$.
Then, it is immediate from Proposition \ref{pr.conv_all} that 
$\lim_{m \to \infty} y(m)_t =y_t$ in $L^{r} ({\cal W},\mu; {\mathbb R}^e)$ for any $r, t$.
Form Proposition \ref{pr.est_ym}, Remark \ref{re.mxi}, and Proposition \ref{pr.conv_all},
$\{  D^n y(m)_t \}_{m=1}^{\infty}$ is Cauchy in 
$L^r ({\cal W},\mu; {\cal H}^{* \otimes n} \otimes {\mathbb R}^e )$ for any $n, r, t$.
Since $D$ is a closed operator, we have 
$
D^n y_t = \lim_{m \to \infty} D^n y(m)_t
$
in $L^r ({\cal W},\mu; {\cal H}^{* \otimes n} \otimes {\mathbb R}^e )$ for any $n, r, t$.
Therefoere, $y_t \in {\mathbb D}_{\infty} ( {\mathbb R}^e) =  
\cap_{n \ge 0} \cap_{2 \le r <\infty} {\mathbb D}_{r,n} ( {\mathbb R}^e) $ for any $t$.
This completes the proof of Theorem \ref{tm.main}.
\toy

\begin{re}
\begin{sloppypar}
The
keys of this proof are the following two facts:
{\rm (i)}~Convergence of $\{ ({\bf w} (m), {\bf b} (m) )\}$ to $({\bf w}, {\bf b} )$ in $L^{\infty -}$
and almost surely.
{\rm (ii)}~Boundedness of 
$\{ \|J(m)\|_{\infty} + \|K(m)\|_{\infty}  \}$  in $L^{\infty -}$, 
i.e. Lemma \ref{lm.unif_exp_int}.
So even when we do not assume {\bf (H)}, we can prove Theorem \ref{tm.main} in the same way
if we can show both {\rm (i)} and {\rm (ii)} by other means.
\end{sloppypar}
\end{re}

%%%%%%%%%%%%%
%\newpage
%%%%%%%%%%%%%%

\section{Proof of Lemma \ref{lm.unif_exp_int} }
\label{sec.pr.lm}

In this section we prove Lemma \ref{lm.unif_exp_int}, following Cass, Litterer, and Lyons \cite{cll}.

For $\alpha >0$ and ${\bf x} \in G\Omega_p ({\mathbb R}^{d})$,
we set $\tau_0 (\alpha) =0$ and 
$$
\tau_{i+1}  (\alpha)
=\inf \{ t \in (\tau_i (\alpha)  ,1]
~|~ 
\omega_{{\bf x}} ( \tau_i (\alpha)  ,t)    
\ge \alpha\} \wedge 1. 
$$
and define 
\begin{equation}
\label{def.locvar2}
N_{\alpha} ({\bf x}) 
=
\sup \{ n \in {\mathbb N} ~|~ \tau_n (\alpha)  <1 \}.
\end{equation} 
Note that
$N_{\alpha} ({\bf x}) $ is non-increasing in $\alpha$ 
and 
$\alpha N_{\alpha} ({\bf x})  \le  \omega_{{\bf x}} (0,1)$
from superadditivity of $\omega_{{\bf x}} $.

$J$ and $K$ can be regarded as solutions of linear RDEs 
driven by a matrix-valued 
rough path ${\bf M} = \int \nabla \sigma ({\bf y}) \la \,\cdot\, , d{\bf x} \ra$. 
For any $\beta >0$, there exists $\alpha >0$ 
such that 
$\omega_{ {\bf M} } (S,T) \le \beta$ if $\omega_{ {\bf x} } (S,T) \le \alpha$ for all $[S,T] \subset [0,1]$.
Growth of the first level paths of solutions of linear RDEs are studied 
in Section 10.7, \cite{fvbk}.
Combining this with Proposition 4.11, \cite{cll} (or the multiplicative property of $J$ and $K$), 
 we have the following deterministic estimate:
for any sufficiently small  $\alpha >0$, there exists a positive constant $C=C_{\alpha}$ such that
\begin{equation}\label{cll.estJ.ineq}
\|  J ({\bf x}) \|_{\infty} +  \|  K ({\bf x}) \|_{\infty} 
\le 
C \exp (  C N_{\alpha} ({\bf x})).
\end{equation}
%
%
%$ {\bf K}^1 ({\bf x})  $ satisfies a similar estimate, too.
%Since $\| {\bf w} (m)\|_{p -var} ~(m=1,2,\ldots)$ are bounded in $L^{\infty -}$and 
%
%
Since $\|J(m)\|_{\infty} +\|K(m)\|_{\infty}  \in  L^{\infty -}$
for each fixed $m$ (as is explained in Section \ref{sec.pwl}),
the problem reduces to showing the following:
for any sufficiently small  $\alpha >0$, there exists $m_0 \in {\mathbb N}$ such that 
$\{ \exp (   N_{\alpha} ({\bf w} (m)))  \}_{ m_0 \le m <\infty}$
are bounded in $L^{\infty -}$.

%Under {\bf (H)},  there exists $p \in (2\rho, 4)$ and $q \in [1,2)$ 
%such that $1/p +1/q >1$
%and ${\cal H}$ is continuously embedded in $C^{q -var}_0 ([0,1], {\mathbb R}^d)$.
%(See Corollary 5.5, \cite{cll}).
%
%As a result, Young translation $\tau_h ({\bf w}) = "{\bf w} +h"$ is well-defined.

%Let $\pi_m$ be the projection $w \mapsto w(m)$and
%
% 
Let $({\cal W}_m ,  {\cal H}_m, \mu_m)$ be the image 
of $({\cal W} ,  {\cal H}, \mu)$ by the projection $w \mapsto w(m)$.
%
%
%$\pi_m$.
%
%
The natural lift will be denoted by ${\cal L}$ below.
We apply Proposition 6.2, \cite{cll} to the new abstract Wiener space.
Then, we have the following:
for any $\hat{w} \in {\cal W}_m$,  $\hat{h} \in {\cal H}_m$, and $\alpha >0$
such that $\omega_{{\cal L} (\hat{w} - \hat{h}) } (0,1) 
:=\sum_{i=1}^{[p]}  \|{\cal L} (\hat{w} - \hat{h})^1  \|^{p/i}_{p/i -var} 
\le \alpha^p $, 
the inequality $\|  \hat{h}\|_{q -var} \ge \alpha N_{ \tilde{\alpha}^p} ({\cal L} (\hat{w}) )$ holds.
Here, 
$\tilde{\alpha} =c_{p,q} \alpha$ for some positive constant 
$c_{p,q} $ which depends only on $p, q$.

Set $B_{\alpha} =\{ {\bf x} \in G\Omega_p ({\mathbb R}^d) 
~|~  \omega_{{\bf x}} (0,1)^{1/p} \le \alpha \}$.
Set also 
$A_{\alpha} =\{ w \in {\cal W} ~|~ {\cal L} (w)={\bf w} \in  B_{\alpha} \}$ 
and 
$A_{\alpha}^{(m)} =\{ w \in {\cal W} ~|~ {\cal L} (w(m))= {\bf w}(m) \in  B_{\alpha} \}$.
By the support theorem for Gaussian rough paths (Theorem 15.60, \cite{fvbk}),
$\mu (A_{\alpha} ) >0$ for any $\alpha >0$.
Since ${\bf w}(m) $ converges to ${\bf w}$ in law, 
there exists a constant $\beta_{\alpha} >0$ such that $\mu (A_{\alpha}^{(m)} ) \ge \beta_{\alpha}$
for sufficiently large $m$.

For any $m$, we have 
$\|h(m) \|_{q -var } \le 3^{1 -1/q} \|h \|_{q -var }  \le 3^{1 -1/q}  C_{{\rm emb}} \|h \|_{{\cal H} }$,
where $C_{{\rm emb}} >0$ 
is the operator norm of the injection.
(See Proposition 5.20, \cite{fvbk} for the first inequality.)
Set $r_n= \alpha 3^{1/q -1} n^{1/q} / C_{{\rm emb}} $.
If $w \in A_{\alpha}^{(m)}$ and $h \in r_n {\cal K}$, where ${\cal K}$ is the unit ball in ${\cal H}$,
then we have
 \[
N_{ \tilde{\alpha}^p} ({\cal L} ( (w+h) (m)) )
=
N_{ \tilde{\alpha}^p} ({\cal L} ( w (m)   +h(m)  ) )
\le 
\frac{\|h(m) \|_{q -var }^q  }{ \alpha^q } \le  n.
\]
This implies that $A_{\alpha}^{(m)} + r_n {\cal K}  \subset
 \{  w \in {\cal W} ~|~  N_{ \tilde{\alpha}^p} ( {\bf w} (m) ) \le n\}$.
Therefore, it follows from Borell's inequality (see Theorem 6.1, \cite{cll}) that
\[
 \{  w \in {\cal W} ~|~  N_{ \tilde{\alpha}^p} ( {\bf w} (m) ) > n\}
\le 1 -\Phi (\gamma_{\alpha} + r_n),
\]
where $\Phi (z) = (2\pi)^{-1/2} \int_{- \infty}^z  e^{-s^2/2 } ds$ 
is the standard normal cumulative ditribution function 
and $\gamma_{\alpha} \in \Phi^{-1} (\beta_{\alpha}) \in {\mathbb R}$.
Here, 
the right hand side is already independent of $m$.
So, we can 
argue in the same way as in Theorem 6.3, \cite{cll}, 
and show
the right hand side is dominated by $C_{\alpha}' \exp (- C_{\alpha}  n^{2/q})$, 
for some positive constants $C_{\alpha}, C_{\alpha}' $ which are independent of $m$.
Hence,  for any $r >0$ and $\alpha>0$,
\[
\sup_{m_0 \le m <\infty} 
 E[ \exp (r N_{ \tilde{\alpha}^p} ( {\bf w} (m)   )  ] < \infty.
\]
Here, $m_0$ depends on $\alpha$, but not on $r$.
Since $\tilde{\alpha}^p$ can be arbitrarily small,   we complete the proof of Lemma \ref{lm.unif_exp_int}.
\toy

%\begin{re}
%\end{re}

%%%%%%%%%%%%%%%%%%%%%%%%%%%%
%%  reference 
%%%%%%%%%%%%%%%%%%%%%%%%%\newpage 

\end{document}